\documentclass{amsart}
\usepackage{graphicx}
\usepackage{amsmath}
\usepackage{amsthm}
\usepackage{amsfonts}
\usepackage{amssymb, amscd}

\usepackage[all]{xy}

\newtheorem{thm}{Theorem}[section]

\newtheorem{lem}[thm]{Lemma}
\newtheorem{proposition}[thm]{Proposition}
\newtheorem{example}[thm]{Example}
\theoremstyle{definition}
\newtheorem{definition}[thm]{Definition}
\theoremstyle{remark}
\newtheorem{remark}[thm]{Remark}
\numberwithin{equation}{section}

\newcommand{\K}{\mathbb K}

\newcommand{\A}{\mathcal{A}}
\newcommand{\G}{\mathcal{G}}

\newcommand{\sll}{\mathfrak{sl}_2(\mathbb{K})}

\begin{document}

\title[  Notes on Formal Deformations of Hom-associative and Hom-Lie Algebras]
{ Notes on Formal Deformations of Hom-associative \\ and Hom-Lie Algebras}%
\author{A. MAKHLOUF and S. SILVESTROV}%
\address{Abdenacer Makhlouf, Universit\'{e} de Haute Alsace,  Laboratoire de Math\'{e}matiques, Informatique et Application,
4, rue des Fr\`{e}res Lumi\`{e}re F-68093 Mulhouse, France}%
\email{Abdenacer.Makhlouf@uha.fr}
\address{Sergei Silvestrov, Centre for Mathematical Sciences,  Lund University, Box
   118, SE-221 00 Lund, Sweden}
\email{sergei.silvestrov@math.lth.se}

\thanks {This work was partially supported by the Crafoord foundation, The Royal Physiographic Society in Lund,
The Swedish Royal Academy of Sciences, The Swedish Foundation of
International Cooperation in Research and High Education (STINT), University of Mulhouse and Lund University  and
the European network Liegrits}
\subjclass[2000]{16S80,16E40,17B37,17B68}
\keywords{Deformation,  Hom-Associative algebra, Hom-Lie
algebra,Hom-Poisson algebra, cohomology}
\date{September 2007}
%
\begin{abstract}
The aim of this paper is to extend to  Hom-algebra structures the
 theory of formal deformations of algebras which was introduced by
Gerstenhaber for associative algebras and extended to Lie algebras
by Nijenhuis-Richardson. We deal with Hom-associative and Hom-Lie
algebras. We construct the first groups of a deformation
cohomology and give several examples of deformations. We provide
families of Hom-Lie algebras deforming Lie algebra $\sll$ and describe as formal deformations the $q$-deformed Witt algebra and Jackson  $\sll$.
\end{abstract}
\maketitle

\section*{Introduction}

Discretization of vector fields via twisted derivations leads to
Hom-Lie and quasi-Hom-Lie structures. This quasi-deformation method
was devised in \cite{HLS,LS1,LS2}. We have introduced in  \cite{MS},
the structure of Hom-associative algebras generalizing associative
algebras to a situation where associativity law is twisted by a linear space
homomorphism. This provided a different way for constructing Hom-Lie
algebras, introduced in \cite{HLS}, by extending the fundamental
construction of Lie algebras from associative algebras via
commutator bracket multiplication. We showed that the commutator
product defined using the multiplication in a Hom-associative
algebra leads naturally to Hom-Lie algebras.  The notion,
constructions and properties of the enveloping algebras of Hom-Lie
algebras are yet to be properly studied in full generality. An
important progress in this direction has been made recently
by D. Yau \cite{Yau:EnvLieAlg}.

The idea to deform algebraic, analytic and geometric structures
within the appropriate category is obviously not new. The first
modern appearance is often attributed to Kodaira and Spencer
and deformations in connection to complex structures on complex manifolds. This
was, however, soon extended and generalized in an
algebraic-homological setting by Gerstenhaber, Grothendieck and
Schlessinger. Nowadays deformation-theoretic ideas penetrate most
aspects of both mathematics and physics and cut to the core of
theoretical and computational problems. In the case of Lie algebras,
quantum deformations (or $q$-deformations) and quantum groups
associated to Lie algebras has been investigated for over twenty
years, still growing richer by the minute. This area began a period
of rapid expansion around 1985 when Drinfel'd  and Jimbo
 independently considered deformations of
$\mathcal{U}(\mathfrak{g})$, the universal enveloping algebra of a
Lie algebra $\mathfrak{g}$, motivated, among other things, by their
applications to the Yang-Baxter equation and quantum inverse
scattering methods.
 The formal deformation theory
was introduced in 1964 by Gerstenhaber \cite{Gerstenhaber64} for
 associative algebras, and in 1967 by Nijenhuis and Richardson
for Lie algebras \cite{Ni-Ri}. In this theory the scalar's field is
extended to the power series ring. The fundamental results of
Gerstenhaber's theory connect deformation theory with the suitable
cohomology groups. There is no general deformation cohomology
theory. Other approaches to study deformation exist
\cite{BjarLaudal, Fialowski86, Fialowski88, Fialowski90, Goze88,
Goze-Remm02,LaudalLNM}, see \cite{Makhlouf-def} for a review.

In this paper we aim to extend the formal deformation theory to
Hom-associative and Hom-Lie algebras and start to construct a
cohomology theory adapted to this deformation theory. In Section 1,
we recall the basic definitions  of Hom-algebra structures and
introduce a definition of Hom-Poisson algebra  which will play an
interesting role in the deformation theory of commutative
Hom-associative algebras. In Section 2, we introduce formal
deformations for Hom-associative and Hom-Lie algebras. Since the
 Hom-Lie algebras enlarge the category of Lie algebras, then one may expect
 non-trivial deformations for rigid associative and Lie algebras. Indeed, we provide Hom-Lie deformations of the classical Lie algebra $\sll$  and  show that the Jackson
$\sll$ is also a  Hom-Lie formal deformation of $\sll$. In Section 3, we give some elements of cohomology of
Hom-associative algebras which are used in Section 4 to describe
deformations of Hom-associative algebras in terms of cohomology.
Section 5 is dedicated to cohomology and deformations of Hom-Lie
algebras. In Section 6, we provide   families of Hom-Lie algebras deforming
 Lie algebra $\sll$ and describe as formal deformations the $q$-deformed Witt algebra.
\section{Hom-algebra structures and Modules}
  A Hom-algebra
structure is a multiplication on a vector space
where the structure is twisted by a linear space homomorphism.
In the following  we summarize the definitions of
Hom-associative, Hom-Leibniz and Hom-Lie
algebraic structures  \cite{MS} generalizing the
well known associative, Leibniz and
Lie-admissible algebras. Also we introduce notion of Hom-Poisson algebras generalizing Poisson algebras to the Hom-Lie context, and define the notion of modules over Hom-associative algebras.
\subsection{Definitions}
Throughout the article we let  $\mathbb{K}$ be an algebraically closed field of characteristic
$0$ and $V$ be a vector space over $\mathbb{K}$.
\begin{definition}
A \emph{Hom-associative algebra} over $V$ is a
triple $( V, \mu, \alpha) $ where  $\mu : V\times
V \rightarrow V$ is a bilinear map and $\alpha:V
\rightarrow V$ is a linear map, satisfying
\begin{equation}\label{Hom-ass}
\mu(\alpha(x), \mu (y, z))= \mu (\mu (x, y), \alpha (z)).
\end{equation}
\end{definition}

In the language of Hopf algebra a
\emph{Hom-associative algebra} over $V$ is a
linear map $\mu : V\otimes V \rightarrow V$ and a
linear map $\alpha$ satisfying
\begin{equation}\label{Hom-ass}
\mu(\alpha(x)\otimes \mu (y\otimes z))= \mu (\mu (x\otimes y)\otimes
\alpha (z)).
\end{equation}

The tensor product of  Hom-associative algebras
$\left( V_1,\mu _1,\alpha_1  \right) $ and
$\left( V_2,\mu_2,\alpha_2 \right) $ is defined
in an obvious way as the Hom-associative algebra
$\left( V_1\otimes V_2,\mu _1 \otimes
\mu_2,\alpha_1 \otimes \alpha_2 \right) $.

A linear map $\phi\ :V\rightarrow V^{\prime }$ is
a morphism
of Hom-associative algebras if%
$$
\mu ^{\prime }\circ (\phi\otimes \phi)=\phi\circ
\mu \quad \text{ and } \qquad \phi\circ
\alpha=\alpha^{\prime }\circ \phi.
$$

In particular, Hom-associative algebras $\left(
V,\mu,\alpha \right) $ and $\left( V,\mu ^{\prime
},\alpha^{\prime }\right) $ are isomorphic if
there exists a
bijective linear map $\phi\ $such that%
$$
\mu =\phi^{-1}\circ \mu ^{\prime }\circ
(\phi\otimes \phi)\qquad \text{ and }\qquad
\alpha= \phi^{-1}\circ \alpha^{\prime }\circ
\phi.
$$

The Hom-associative algebra is called unital if
it admits a unity, i.e. an element $e\in V$ such
that $\mu(x,e)=\mu(e,x)=x$ for all $x\in V$. The
unity may be also expressed by a linear map
$\eta:\K \rightarrow V$ defined by $\eta(c)= c e$
for all $c\in \K$.
 Let $\left( V,\mu ,\alpha,\eta \right) $ and $\left( V^{\prime },\mu
^{\prime },\alpha^{\prime },\eta ^{\prime
}\right) $ be two unital Hom-associative
algebras. Then morphisms of unital
Hom-associative algebras are required also to
preserve the unital structure, i.e. satisfy
$f\circ \eta =\eta ^{\prime }, $ and in the
definition of isomorphism of unital
Hom-associative algebras it is also required
that $\eta =f^{-1}\circ \eta ^{\prime }.$

The Hom-Lie algebras were introduced by Hartwig,
Larsson and Silvestrov in \cite{HLS} motivated
initially by examples of deformed Lie algebras
coming from twisted discretizations of vector
fields.

\begin{definition}
A \emph{Hom-Lie algebra} is a triple $(V, [\cdot, \cdot], \alpha)$
consisting of
 a linear space $V$, bilinear map $[\cdot, \cdot]: V\times V \rightarrow V$ and
 a linear space homomorphism $\alpha: V \rightarrow V$
 satisfying
$$\begin{array}{c} [x,y]=-[y,x] \quad {\text{(skew-symmetry)}} \\{}
\circlearrowleft_{x,y,z}{[\alpha(x),[y,z]]}=0 \quad
{\text{(Hom-Jacobi condition)}}
\end{array}$$
for all $x, y, z$ from $V$, where $\circlearrowleft_{x,y,z}$ denotes
summation over the cyclic permutation on $x,y,z$.
\end{definition}
In a similar way we have the following definition of Hom-Leibniz
algebra.
\begin{definition}
A \emph{Hom-Leibniz algebra} is a triple $(V, [\cdot, \cdot],
\alpha)$ consisting of a linear space $V$, bilinear map $[\cdot,
\cdot]: V\times V \rightarrow V$ and a homomorphism $\alpha: V
\rightarrow V$  satisfying
\begin{equation} \label{Leibnizalgident}
 [[x,y],\alpha(z)]=[[x,z],\alpha (y)]+
 [\alpha(x),[y,z]].
\end{equation}
\end{definition}
Note that if a Hom-Leibniz algebra is skewsymmetric then it is a
Hom-Lie algebra.

We introduce in the following the definition of Hom-Poisson
structure which involves naturally in the deformation theory.

\begin{definition}
A \emph{Hom-Poisson algebra} is a quadruple $(V,\mu, \{\cdot,
\cdot\}, \alpha)$ consisting of
 a linear space $V$, bilinear maps $\mu: V\times V \rightarrow V$ and
  $\{\cdot, \cdot\}: V\times V \rightarrow V$, and
 a linear space homomorphism $\alpha: V \rightarrow V$
 satisfying
 \begin{enumerate}
\item $(V,\mu, \alpha)$ is a commutative Hom-associative algebra,
\item $(V, \{\cdot,
\cdot\}, \alpha)$ is a Hom-Lie algebra,
\item
for all $x, y, z$ in $V$,
\begin{equation}\label{CompatibiltyPoisson}
\{\alpha (x) , \mu (y,z)\}=\mu (\alpha (y),
\{x,z\})+ \mu (\alpha (z), \{x,y\}).
\end{equation}
\end{enumerate}
\end{definition}
The condition \eqref{CompatibiltyPoisson}
expresses the compatibility between the
multiplication and the Poisson bracket. It can be
 reformulated equivalently as
\begin{equation}\label{CompatibiltyPoissonLeibform}
\{\mu(x,y),\alpha (z) \}=\mu (\{x,z\},\alpha
(y))+\mu (\alpha (x), \{y,z\})
\end{equation}
for all $x, y, z$ in $V$. Note that in this form
it means that $ad_z (\cdot) = \{\cdot,z\}$ is a sort  generalization of
derivation of associative  algebra defined by  $\mu$, and
also it resembles the identity
\eqref{Leibnizalgident} in the definition for
Leibniz algebra.

We also recall in the following the structure
of module over
Hom-associative algebras.
\begin{definition}
Let $\A=(V,\mu,\alpha)$ be a Hom-associative
$\K$-algebra. An $\A$-module (left) is a triple
$(M,f,\gamma)$ where $M$ is $\K$-vector space and
$f,\gamma$ are  $\K$-linear maps, $f:  M
\rightarrow M$ and $\gamma : V \otimes M
\rightarrow M$, such that the following diagram
commute:

$$
\begin{array}{ccc}
V\otimes V\otimes M & \stackrel{\mu \otimes f}{\longrightarrow } &
V\otimes
M \\
\ \quad \left\downarrow ^{\alpha \otimes \gamma }\right. &  & \quad
\left\downarrow
^\gamma \right. \\
V\otimes M & \stackrel{\gamma }{\longrightarrow } & M
\end{array}
$$
\end{definition}
\begin{remark}
A Hom-associative $\K$-algebra $\A=(V,\mu,\alpha)$ is  a left
$\A$-module with $M=V$, $f=\alpha$ and $\gamma =\mu$.
\end{remark}

\section{Formal deformations of  Hom-associative algebras and Hom-Lie algebras}

In this section we extend to Hom-algebra
structures the formal deformation theory
introduced  by Gerstenhaber for associative
algebras \cite{Gerstenhaber64}, and by Nijenhuis
and Richardson for Lie algebras \cite{Ni-Ri}.
More precisely, we define the concept of
deformation for Hom-associative algebras and
Hom-Lie algebras and define the suitable $1^{st}$
and $2^{nd}$ cohomology groups adapted to formal
deformation.

Let $V$ be a $\K$-vector space and ${\A }_{0}=(V,\mu _{0},\alpha
_0)$ be a  Hom-associative algebra   and $L_0=(V, [\cdot, \cdot]_0,
\alpha_0)$  be a Hom-Lie algebra. Let $\K [[t]]$ be the power series
ring in one variable $t$ and coefficients in $\K $ and $V[[t]]$ be
the set of formal power series whose coefficients are elements of $V$, ($V[[t]]$ is obtained   by extending the coefficients domain of $V$ from $\K $
to $\K [[ t]]$). Then $V[[t]]$ is a $\K[[t]]$-module. When $V$ is
finite-dimensional, we have $ V[[t]]=V\otimes _{\K }\K[[t]]$. Note that $V$ is a submodule of $V[[t]]$. Given a
 $\K$-bilinear map $f :V\times V  \rightarrow V$, it admits naturally an extension to a $\K[[t]]$-bilinear map $f :V[[t]]\times V[[t]]  \rightarrow V[[t]]$, that is,  if  $x=\sum_{i\geq0}{a_i t^i}$ and $y=\sum_{j\geq0}{b_j t^j}$ then
$f(x,y)=\sum_{i\geq0,j\geq0}{t^{i+j}f (a_i,b_j)}$. The same holds for  linear maps.

\begin{definition}
Let $V$ be a $\K$-vector space and ${\A }_{0}=(V,\mu _{0},\alpha_0)$ be a Hom-associative  algebra.
A \emph{formal Hom-associative deformation} of ${\A }_{0}$  is given by the
$\K[[t]]$-bilinear and $\K[[t]]$-linear maps
$\mu_{t} :V[[t]]\times V[[t]]  \rightarrow V[[t]]$ and
$\alpha_{t}:V[[t]] \rightarrow  V[[t]]$ of the form
$$\mu_{t} =\sum_{i\geq 0}\mu_{i}t^{i},\quad\text{and}\quad \alpha_{t} =\sum_{i\geq 0}\alpha_{i}t^{i}
$$
where each $\mu_{i}$ is a $\K$-bilinear map $\mu_{i}:  V\times V
\rightarrow  V$ (extended to be $\K[[t]]$-bilinear) and each
$\alpha_{i}$ is a $\K$-linear map $\alpha_{i}:  V  \rightarrow  V$
(extended to be $\K[[t]]$-linear), such that holds for $x, y,z\in
 V$ the following
formal Hom-associativity condition:
\begin{equation}
\label{equ1} \mu_{t}(\alpha_t(x), \mu_{t}(y,
z))-\mu_{t}(\mu_{t}(x, y), \alpha_t(z)))=0.
\end{equation}
\end{definition}

\begin{definition}

Let $V$ be a $\K$-vector space and  $L_0=(V, [~, ~]_0,
\alpha_0)$ be a Hom-Lie algebra. A \emph{formal Hom-Lie deformation} of $L_{0}$  is given by the $\K[[t]]$-bilinear and $\K[[t]]$-linear maps $
[~, ~]_{t}:  V[[t]] \times V[[t]]  \rightarrow
V[[t]], \ \ \alpha_{t}: V[[t]] \rightarrow  V[[t]]
$  of the form
$$[~,~]_{t} =\sum_{i\geq 0}[~,~]_{i}t^{i},\quad\text{and}\quad \alpha_{t} =\sum_{i\geq 0}\alpha_{i}t^{i}
$$
where each $[~,~]_{i}$ is a $\K$-bilinear map $[~,~]_{i}:  V\times V  \rightarrow  V$ (extended to be $\K[[t]]$-bilinear)  and each $\alpha_{i}$ is a $\K$-linear map $\alpha_{i}:  V  \rightarrow  V$ (extended to be $\K[[t]]$-linear), and  satisfying for $x, y,z\in
 V$ the following conditions:
\begin{equation}
\label{equa2}
\begin{array}{c}
[x,y]_{t}=-[y,x]_{t} \quad
{\text{(skew-symmetry)}} \\{}
\circlearrowleft_{x,y,z}{[[\alpha_t(x),[y,z]_t]_t}=0.
\quad {\text{(Hom-Jacobi condition)}}
\end{array}
\end{equation}
\end{definition}

\begin{remark}
The skew-symmetry of $[\cdot,\cdot]_{t}$ is
equivalent to the skew-symmetry of all
$[\cdot,\cdot]_{i}$ for $i\in \mathbb{Z}_{\geq
0}$.
\end{remark}

We call the condition \eqref{equ1} (respectively
\eqref{equa2})  deformation equation of
Hom-associative (respectively Hom-Lie) algebra.

\begin{example}[\textbf{Jackson $\sll$}]

In this example, we will consider the Hom-Lie algebra Jackson $\sll$
which is a Hom-Lie deformation of the classical Lie algebra  $\sll$ defined by $[
h,f] = -2f,\   [ h,e] = 2e, \  [ e,f] =h$. The Jackson $\sll$ is
related to Jackson derivations.  As linear space, it is generated by
$e,f,h$ with the brackets defined by
$$[ h,f]_t = -2f-2tf,\ \  [ h,e]_t = 2e, \ \ [ e,f]_t =h+\frac{t}{2}h.$$

The linear map $\alpha_t$ is defined by
$$\alpha_t(e)=\frac{2+t}{2(1+t)}e=e+\sum_{k=0}^{\infty}{\frac{(-1)^k}{2}t^k\ } e, \quad
\alpha_t(h)= h, \quad \alpha_t(f)=f+\frac{t}{2} \, f.$$
The Hom-Jacobi identity is proved as follows. It is enough to
consider it on $e$, $f$ and $h$:
\begin{align*}
[ \alpha_t(e),&[ f,h]_t]_t+[ \alpha_t(f),[ h,e]_t]_t+[
\alpha_t(h),[ e,f]_t]_t=\\
&=(2+t)[ e,f]_t+(2+t)[ f,e]_t+\frac{(2+t)}{2}[ h,h]_t =0,
\end{align*}
where we have used that $[~,~]_t$ is skew-symmetric.

In this case
$$[ h,f]_1 = -2f,\ \  [ h,e]_1 = 0, \ \ [ e,f]_1 =\frac{1}{2}h,$$
$$\alpha_1(e)=-\frac{1}{2}\  e, \quad
\alpha_1(h)=0 , \quad \alpha_1(f)=\frac{1}{2} \, f.$$ And for
$k\geq2$, one has
$$[ h,f]_k = 0,\ \  [ h,e]_k = 0, \ \ [ e,f]_k =0,$$

$$\alpha_k(e)=\frac{(-1)^k}{2}\  e, \quad
\alpha_k(h)= 0, \quad \alpha_k(f)=0.$$ Thus Jackson $\sll$ algebra
is a Hom-Lie algebra deformation of $\sll$. Indeed,
$$[ h,f]_0 = -2f,\ \  [ h,e]_0 = 2e, \ \ [ e,f]_0 =h,$$
$$\alpha_0(e)=e, \quad
\alpha_0(h)= h, \quad \alpha_0(f)=f.$$
\end{example}
\subsection{Deformation equation of Hom-associative algebras}

In this section, we study the equation (\ref{equ1}) and  thus
characterize the deformations of  Hom-associative algebras. The
equation  may be written

\begin{equation}
\sum_{i\geq 0}\sum_{j\geq 0}\sum_{k\geq 0}(\mu_{i}(\alpha_{k}(x),
\mu_{j}(y,z))-\mu_{i}(\mu_{j}(x, y),\alpha_k (z)))t^{i+j+k}=0
\end{equation}

\begin{definition}
We call $\alpha$-associator the map
\begin{equation}\label{associator}
Hom(V^{\times 2},V)\times Hom(V^{\times 2},V)  \longrightarrow
Hom(V^{\times 3},V), \quad
(\mu_{i},\mu_{j})  \longmapsto \mu_{i}\circ_{\alpha} \mu_{j},%
\end{equation}
defined for all $x,y,z\in V$ by
\begin{equation}
\mu_{i}\circ_\alpha \mu_{j}(x, y, z)
=\mu_{i}(\alpha(x),
\mu_{j}(y,z))-\mu_{i}(\mu_{j}(x, y),\alpha (z)).
\end{equation}
\end{definition}

By using  $\alpha_k$-associators, the deformation equation may be
written as follows
\begin{equation}
\sum_{i\geq 0}\sum_{j\geq 0} \sum_{k\geq
0}(\mu_{i}\circ_{\alpha_k} \mu_{j})t^{i+j+k}=0
\quad\text{or}\quad \sum_ {s\geq 0}t^{s}\sum_ {k=
0}^{s}\sum_ {i=0}^{s-k}\mu_{i}\circ_{\alpha_k}
\mu_{s-k-i}=0.
\end{equation}
This equation is equivalent to the following
infinite system:
\begin{equation}\label{DefEqu}
\sum_ {k= 0}^{s}\sum_ {i=0}^{s-k}\mu_{i}\circ_{\alpha_k}
\mu_{s-k-i}=0,\quad s=0,1,\ldots
\end{equation}
In particular,
\begin{itemize}
\itemsep-3pt
\item
for $s=0$, \ $\mu_{0}\circ_{\alpha_0} \mu_{0}=0,$
this corresponds to the Hom-associativity of
$\A_{0}$;
\item
for $s=1$, \ $\mu_{1}\circ_{\alpha_0}
\mu_{0}+\mu_{0}\circ_{\alpha_1}
\mu_{0}+\mu_{0}\circ_{\alpha_0} \mu_{1}=0$;
\item
for $s=2$, \ $\mu_{2}\circ_{\alpha_0}
\mu_{0}+\mu_{0}\circ_{\alpha_0} \mu_{2}+\mu_{1}\circ_{\alpha_0}
\mu_{1}+\mu_{0}\circ_{\alpha_1} \mu_{1}+\mu_{1}\circ_{\alpha_1}
\mu_{0}+\mu_{0}\circ_{\alpha_2} \mu_{0}=0$.
\end{itemize}

\subsection{Equivalent and trivial deformations}

In this section, we characterize the equivalent and trivial
deformations of Hom-associative algebras.

\begin{definition}
Let  $\A_0=(V,\mu_0,\alpha_0)$ be a Hom-associative algebra. Given
two deformations of $\A_0$, $\A_t=(V,\mu_t,\alpha_t)$ and
$\A'_t=(V,\mu'_t,\alpha'_t)$  where
 $\mu_{t}=\sum_{i\geq 0}\mu_{i}t^{i}$, $\mu'_{t}=\sum_{i\geq 0}\mu'_{i}t^{i}$,
  $\alpha_{t}=\sum_{i\geq 0}\alpha_{i}t^{i}$ and $\alpha'_{t}=\sum_{i\geq
  0}\alpha'_{i}t^{i}$, with
 $\mu_{0}=\mu'_{0}$ and $\alpha_{0}=\alpha'_{0}$.

 We say
that they are \emph{equivalent} if there is a formal isomorphism
$\Phi_{t}: V\rightarrow V[[t]]$ which is a $\K[[t]]$-linear map that
may be written in the form $ \Phi_{t}=\sum_{i\geq 0}\Phi _{i}t^{i}
=Id+\Phi _{1}t+\Phi_{2}t^{2}+\ldots$  where $\Phi_{i}\in End_{\K
}(V)$ and $\Phi_{0}=Id $ such that
\begin{equation}
\label{equIso} \Phi_{t}\circ
\mu_{t}=\mu'_{t}\circ(\Phi _{t}\times \Phi _{t}
)\quad \text{and} \quad \alpha'_t \circ \Phi _{t}
= \Phi _{t}\circ \alpha_t.
\end{equation}
A deformation $\A_{t}$ of $\A_{0}$ is said to be \emph{trivial} if
and only if $\A_{t}$ is equivalent to $\A_{0}$, viewed  as an
algebra over $V[[t]]$.
\end{definition}

The condition (\ref{equIso}) may be written
\begin{equation}\label{isom1}
\Phi _{t}(\mu_{t}(x, y)) =\mu'_{t}(\Phi _{t}(x)), \Phi
_{t}(y)),\quad\forall x,y\in V
\end{equation}
and
\begin{equation}\label{isom2}
\Phi _{t}(\alpha_{t}(x)) =\alpha'_{t}(\Phi
_{t}(x)),\quad\forall x\in V.
\end{equation}

The equation \ref{isom1} is equivalent to
\begin{equation}
\sum_{i\geq 0}\Phi _{i}\left(\sum_{j\geq 0}\mu_{j} (x,
y)t^{j}\right)t^{i} =\sum_{i\geq 0}\mu'_{i} \left( \sum_{j\geq
0}\Phi _{j}(x)t^{j},\sum_{k\geq 0}\Phi_{k}(y)t^{k}
 \right) t^{i}
\end{equation}
or
\begin{equation*}
\sum_{i,j\geq 0}\Phi _{i}(\mu_{j}(x, y))t^{i+j} =
\sum_{i,j,k\geq 0}\mu'_{i}(\Phi _{j}(x),
\Phi_{k}(y))t^{i+j+k}.
\end{equation*}
By identification of  coefficients, one obtains
that the constant coefficients are identical,
i.e.
\begin{equation*}
\mu_{0}=\mu'_{0} \quad\text{because}\quad \Phi
_{0} =Id.
\end{equation*}
For  coefficients of $t$ one has
\begin{equation}
\Phi _{0}(\mu_{1}(x, y))+\Phi _{1}(\mu_{0}(x, y))
= \mu'_{1}(\Phi _{0}(x), \Phi _{0}(y)) +
\mu'_{0}(\Phi _{1}(x), \Phi _{0}(y)) +
\mu'_{0}(\Phi _{0}(x), \Phi _{1}(y)).
\end{equation}
Since $\Phi _{0}=Id$, it follows that
\begin{equation}
\mu_{1}(x, y)+\Phi _{1}(\mu_{0}(x, y) =
\mu'_{1}(x, y)+\mu_{0}(\Phi _{1}(x), y) +
\mu_{0}(x, \Phi _{1}(y)).
\end{equation}
Consequently,
\begin{equation}\label{equiv1}
\mu'_{1}(x, y) = \mu_{1}(x, y)+\Phi_{1}(\mu_{0}(x, y))-\mu_{0}(\Phi
_{1}(x), y) - \mu_{0}(x, \Phi _{1}(y)).
\end{equation}

The condition on homomorphisms  \eqref{isom2} is
equivalent for $x\in V$ to
\begin{equation}
\sum_{i,j\geq 0}{\Phi _{i}(\alpha_j (x))
t^{i+j}}=\sum_{i,j\geq 0}{ \alpha'_{i}(\Phi_j
(x)) t^{i+j}}.
\end{equation}
The condition  implies that $\alpha_0=\alpha'_0
\mod t$ and  that
\begin{equation}
\Phi _{1}\circ \alpha_0+\Phi _{0}\circ
\alpha_1=\alpha'_0\circ \Phi
_{1}+\alpha'_1\circ\Phi _{0} \mod t^2.
\end{equation}
Then
\begin{equation}\label{equiv2}
\alpha'_1=\alpha_1+\Phi _{1}\circ
\alpha_0-\alpha_0 \circ \Phi _{1}.
\end{equation}
Then, The first and second order conditions of the equivalence
between two deformations of a Hom-associative algebra are given by
\eqref{equiv1} and \eqref{equiv2}.
\section{First and Second cohomology groups of a Hom-associative algebras}

We introduce in the following  certain elements of the cohomology of
Hom-associative algebras which fits with the deformation theory.

\subsection{First and second coboundary operators}

 Let ${\A }=(V,\mu,\alpha)$ be a Hom-associative algebra on
a $\K$-vector space $V.$

The set of $p$-cochains on $V$ is the set of $p$-linear maps
\begin{equation*}
\mathcal{C}^{p}({\A,\A})=\{\varphi:V^{\times p}
=\underset{p\text{ times}}{\underbrace{V\times
V\times ...\times V}} \longrightarrow V\}.
\end{equation*}
\begin{definition}We set for a morphism $\tau\in Hom(V,V)=\mathcal{C}^{1}({\A,\A})$
$$\rho_{\A}^1 : \mathcal{C}^{1}({\A,\A})\longrightarrow
\mathcal{C}^{1}({\A,\A}) \ \ \  \ \tau \rightarrow \rho_{\A}^1
\tau:= \tau\circ\alpha-\alpha\circ \tau,$$ and
$$\rho_{\A}^2 : \mathcal{C}^{2}({\A,\A})
\longrightarrow \mathcal{C}^{3}({\A,\A})\ \ \  \ \tau \rightarrow
\rho_{\A}^2 \tau :=\mu\circ_\tau\mu.$$
\end{definition}
Now, we define $1$-Hom-cochains and
$2$-Hom-cochains of $\A$.
\begin{definition}
A $1$-\emph{Hom-cochain} of $\A$ is a map $f$,
where $f\in \mathcal{C}^{1}({\A,\A})$  satisfies
\begin{equation}
\rho_{\A}^1 f=0.
\end{equation}
We denote by $Hom\mathcal{C}^{1}({\A,\A})$ the set of all
$1$-\emph{Hom-cochains} of $\A$.
\end{definition}

\begin{definition}
A $2$-\emph{Hom-cochain} is a pair $(\varphi, \tau)$,
where
$\varphi\in \mathcal{C}^{2}({\A,\A})$ and $\tau$
is a linear map
such that
\begin{equation}
\rho_{\A}^2 {\tau}=0.
\end{equation}
We denote by $Hom\mathcal{C}^{2}({\A,\A})$ the
set of all
$2$-\emph{Hom-cochains} of $\A$.
\end{definition}

The $1$-coboundary and $2$-coboundary operators
for Hom-associative algebras are defined as
follows.

\begin{definition}
We call  $1$-\emph{coboundary} operator
of Hom-associative algebra
${\A}$ the map
\begin{equation*}
\delta_{Hom} ^{1}:  \mathcal{C} ^{1}({\A,\A})
\longrightarrow
\mathcal{C}^{2}({\A ,\A }), \quad f  \longmapsto
\delta_{Hom}^{1}f
\end{equation*}
defined by
\begin{equation*}
\delta_{Hom} ^{1}f(x, y) =f(\mu(x,y))-\mu(f(x),
y)- \mu(x, f(y))).
\end{equation*}

\end{definition}

\begin{definition}
We call a $2$-\emph{coboundary} operator
of Hom-associative algebra
${\A}$ the map
\begin{equation*}
\delta_{Hom}^{2}:  \mathcal{C}^{2}({\A ,\A }) \longrightarrow
\mathcal{C}^{3}({\A ,\A }), \quad
 \varphi  \longmapsto  \delta_{Hom} ^{2}\varphi
\end{equation*}
defined by
\begin{align*}
\delta_{Hom}^{2}\varphi (x, y, z) &= \varphi(\alpha(x),
\mu(y,z))-\varphi(\mu(x, y),\alpha (z))\\
&+\mu(\alpha(x), \varphi(y,z))- \mu(\varphi(x,
y),\alpha (z)).
\end{align*}
\end{definition}

\begin{remark}
The operator $\delta_{Hom}^{2}$ can also be
defined using the $\alpha$-associator
\eqref{associator} by
\begin{equation*}
\delta_{Hom}^{2}\varphi =\varphi \circ_\alpha \mu
+ \mu \circ_\alpha \varphi.
\end{equation*}
\end{remark}

The cohomology spaces relative to these
coboundary operators are:

\begin{definition}
The space of $1$-cohomology classes of ${\A}$ is
\begin{equation*}
\mathit{H_{Hom}^{1}}({\A,\A})=\{f\in Hom\mathcal{C}^{1}({\A ,\A
})\,:\,\delta_{Hom}^{1}f=0\ \}=\{f\in \mathcal{C}^{1}({\A ,\A
})\,:\,\delta_{Hom}^{1}f=0\ \text{ and } \rho_{\A} ^1 f=0\}.
\end{equation*}
The space of $2$-coboundaries of ${\A }$ is
\begin{equation*}
\mathit{B_{Hom}^{2}}({\A ,\A }) =\{(\varphi,\tau)\in
Hom\mathcal{C}^{2}({\A ,\A }) \ :\, \varphi = \delta_{Hom} ^{1}f,\
f\in Hom\mathcal{C} ^{1}({\A },{\A })\ \text{ and } \rho_{\A}^2 \tau
 =0\}.
\end{equation*}
The space of $2$-cocycles of ${\A}$ is
\begin{equation*}
\mathit{Z_{Hom}^{2}}({\A ,\A }) =\{(\varphi,\tau)\in
Hom\mathcal{C}^{2}({\A ,\A })\,\,:\,\delta_{Hom} ^{2}\varphi=0 \
\text{ and } \rho_{\A}^2 \tau  =0 \}.
\end{equation*}
\end{definition}
\begin{proposition}
$\delta_{Hom} ^{2}( \delta_{Hom} ^{1})=0.$
\end{proposition}
\begin{proof}
Let
\begin{equation*}
\delta_{Hom} ^{1}f(x, y) =f(\mu(x,y))-\mu(f(x),
y)- \mu(x, f(y))).
\end{equation*}
Then
\begin{align*}
\delta_{Hom}^{2}(\delta_{Hom} ^{1}f )(x, y, z)
&= \delta_{Hom}
^{1}f(\alpha(x),
\mu(y,z))-\delta_{Hom} ^{1}f(\mu(x, y),\alpha (z))\\
&+\mu(\alpha(x), \delta_{Hom} ^{1}f(y,z))-\mu(\delta_{Hom} ^{1}f(x,
y),\alpha (z))\\ &=f(\mu(\alpha(x), \mu(y,z)))-\mu(f(\alpha(x)),
\mu(y,z))-\mu(\alpha(x), f(\mu(y,z)))\\& -f(\mu(\mu(x,y),\alpha(z)
))+\mu(f(\mu(x,y)),\alpha(z) )+\mu(\mu(x,y), f(\alpha(z)))\\&+
\mu(\alpha(x),f( \mu(y,z)))-\mu(\alpha(x)),
\mu(f(y),z)-\mu(\alpha(x), \mu(y,f(z)))\\&
-\mu(f(\mu(x,y)),\alpha(z) )+\mu(\mu(f(x),y)),\alpha(z)
)+\mu(\mu(x,f(y)), \alpha(z))\\&=0,
\end{align*}
because $\alpha$ and $f$ commute and the
multiplication $\mu$ is
Hom-associative.
\end{proof}
\begin{remark}
One has $\mathit{B_{Hom}^{2}}({\A ,\A }){\subset
}\mathit{Z_{Hom}^{2}}({\A ,\A })$,
because $\delta_{Hom} ^{2}\circ
\delta_{Hom} ^{1}=0$. Note also
that $\mathit{H_{Hom}^{1}}({\A ,\A
})$ corresponds to the derivations space
of a Hom-associative
algebra ${\A }$.
\end{remark}

\begin{definition}
We call the $2^{th}$ cohomology group
of the Hom-associative algebra
${\A }$,  the quotient
\begin{equation*}
\mathit{H_{Hom}^{2}}(\A ,\A )
=\frac{\mathit{Z_{Hom}^{2}}(\A
,\A)}{\mathit{B_{Hom}^{2}}(\A ,\A )}.
\end{equation*}
\end{definition}

\begin{remark}
The cohomology class of an element $(\varphi,\tau)$ is given by the set of elements of the form $(\psi,\tau)$ such that $\psi=\delta_{Hom} ^{2}f$ where $f$ is a $1$-Hom-cohain, that is $f\in\mathcal{C} ^{1}({\A },{\A })$ and $\rho_\A ^{1}f=0$.
\end{remark}
\section{Cohomological approach of Hom-associative algebra deformations}
Let $\A_{t}=(V,\mu_{t},\alpha_{t})$  be a deformation of a
Hom-associative algebra ${\A }_{0}=(V,\mu_{0},\alpha_{0})$ where $
\mu_{t}(x, y) =\sum_{i\geq 0}\mu_{i}(x, y)t^{i}$ and $ \alpha_{t}(x)
=\sum_{i\geq 0}\alpha_{i}(x)t^{i}$. We characterize under the
assumption $\rho_{\A}^2{\alpha_{i}} =\mu_0\circ_{\alpha_i}\mu_0=0$
for $i\geq 1 $ the deformations of ${\A }_{0}$ in terms of
cohomology.

By using the definition of $2$-coboundaries and by gathering some
terms, the deformation equation (\ref{DefEqu}) may be written
$$\delta_{Hom} ^{2}\mu_{1}=0
$$
and
\begin{equation}\label{obstruct}
\delta_{Hom} ^{2}\mu_{s}=-\sum_{k=1}^{s-1} {
\sum_{p=0}^{k}\mu_{s-k}\circ_{\alpha_{p}}
\mu_{k-p}}-\sum_{k=1}^{s-k} {\mu_{0}\circ_{\alpha_k}
\mu_{s-k}},\quad s=2,3\ldots .
\end{equation}
Consequently, the following Lemma holds.
\begin{lem}
The pair $(\mu_{1},\alpha_1)$  of the deformation $\A_{t}$ is a
$2$-Hom-cocycle of the cohomology of the Hom-associative algebra
${\A }_{0}$.
\end{lem}

\begin{definition}
Let ${\A }_{0}=(V,\mu_{0},\alpha_0)$ be a Hom-associative algebra
and $(\mu_{1},\alpha_1)$ be an element of
$\mathit{Z}^{2}_{Hom}(\A_{0},\A_{0})$. The 2-Hom-cocycle
$(\mu_{1},\alpha_1)$ is said \emph{integrable} if there exists a
 a pair $(\mu_{t},\alpha_{t})$ such that
$\mu_{t}=\sum_{i\geq 0}\mu_{i}t^{i}$ and $\alpha_{t}=\sum_{i\geq
0}\alpha_{i}t^{i}$ defining a deformation
$\A_{t}=(V,\mu_{t},\alpha_{t})$ of ${\A }_{0}$.
\end{definition}

\begin{proposition}
Let $\A_{t}=(V,\mu_{t},\alpha_{t})$ be a deformation of a
Hom-associative algebra $\A_0$. Let $\mu_{t}=\sum_{i\geq
0}\mu_{i}t^{i}$, $\alpha_{t}=\sum_{i\geq 0}\alpha_{i}t^{i}$ and
$(\mu_{1},\alpha _{1})$ be an element of $\mathit{Z_{Hom}^{2}}(\A_0
,\A_0)$. The integrability of $(\mu_{1},\alpha _{1})$ depends only
on its cohomology class.
\end{proposition}

\begin{proof}
We saw in Section 2 that if two deformations
$\A_{t}=(V,\mu_{t},\alpha_{t})$ and
$\A'_{t}=(V,\mu'_{t},\alpha'_{t})$ are equivalent then

\begin{equation*}
\mu'_{1}(x, y) = \mu_{1}(x, y)+\Phi_{1}(\mu_{0}(x, y))-\mu_{0}(\Phi
_{1}(x), y) - \mu_{0}(x, \Phi _{1}(y))
\end{equation*}
and
\begin{equation*}
\alpha'_1=\alpha_1+\Phi _{1}\circ \alpha_0-
\alpha_0\circ \Phi _{1}.
\end{equation*}
With $\alpha'_1=\alpha_1$, these conditions means that
 $$\mu'_{1}=\mu_{1}+\delta^1_{Hom} \Phi _{1}\quad \text{ and } \quad \rho_\A ^1 (\Phi)=0.$$
Therefore the two elements are cohomologous.

Thus,
\begin{align*}
\delta^2_{Hom} \mu_{1}=0 &\quad\Longrightarrow\quad \delta^2_{Hom}
\mu'_1=\delta^2_{Hom} (\mu_{1}+\delta^1_{Hom} \Phi
_{1})=\delta^2_{Hom} \mu_{1}+\delta^2_{Hom}
(\delta^1_{Hom} \Phi _{1})=0;\\
\mu_{1}=\delta^1_{Hom} g &\quad\Longrightarrow\quad
\mu'_{1}=\delta^1_{Hom} g-\delta^1_{Hom} \Phi_{1}=\delta^1_{Hom}
(g-\Phi _{1});
\end{align*}
which ends the proof.
\end{proof}

\begin{proposition}
Let ${\A }_{0}=(V,\mu_{0},\alpha_0)$ be a Hom-associative algebra.
There is, over $\K [[t]]/t^2$, a one-to-one correspondence between
the elements of $\mathit{H_{Hom}^{2}}({\A }_{0}{ ,\A}_{0})$ and the
infinitesimal deformation of ${\A }_{0}$ defined by
\begin{equation}
\mu_{t}(x, y) =\mu_{0}(x, y)+\mu_{1}(x, y)t, \quad \text{ and }\quad
\alpha_{t}(x) =\alpha_{0}(x)+\alpha_{1}(x)t \quad \forall x,y\in V
\end{equation}
with $\rho_\A ^2 (\alpha_1)=\mu_0\circ_{\alpha_1}\mu_0=0$
\end{proposition}
\begin{proof} The deformation equation is equivalent to
$$\delta^2_{Hom} \mu_{1}+\rho_\A ^2 (\alpha_1)=0.$$
Since $\rho_\A ^2 (\alpha_1)=0$, then the previous equation is
equivalent to $(\mu_{1},\alpha _{1})\in\mathit{Z_{Hom}^{2}}(\A_0
,\A_0)$.
\end{proof}

\subsection{Poisson algebra}
We consider now a commutative Hom-associative algebra and show that
first order deformation induces a Hom-Poisson structure (Definition
1.4). More generally, the following lemma shows that  any
skewsymmetric  2-Hom-cocycle of  a commutative algebra  satisfies
the compatibility condition \ref{CompatibiltyPoisson} with the
multiplication of the Hom-associative algebra.
\begin{lem}
Let  ${\A }=(V,\mu,\alpha)$ be a commutative Hom-associative algebra
and $\varphi$ be a skewsymmetric $2$-cochain such that
$\delta_{Hom}^{2}\varphi=0$. Then for $x,y,z\in V$
\begin{equation}
\varphi(\alpha (x),\mu(y,z))=\mu(\alpha(y),\varphi(x,z))+
\mu(\alpha(z),\varphi(x,y)).
\end{equation}
\end{lem}
\begin{proof}
The condition $\delta_{Hom}^{2}\varphi=0$ is equivalent to
\begin{align*}
 \varphi(\alpha(x),
\mu(y,z))-\varphi(\mu(x, y),\alpha (z))+\mu(\alpha(x),
\varphi(y,z))-\mu(\varphi(x, y),\alpha (z))=0
\end{align*}
Then one has
\begin{equation}\label{eq1}
\varphi(\alpha(x), \mu(y,z))=\varphi(\mu(x, y),\alpha
(z))-\mu(\alpha(x), \varphi(y,z))+\mu(\varphi(x, y),\alpha (z))
\end{equation}
\begin{equation}\label{eq2}
\varphi(\alpha(x), \mu(z,y))=\varphi(\mu(x, z),\alpha
(y))-\mu(\alpha(x), \varphi(z,y))+\mu(\varphi(x, z),\alpha (y))
\end{equation}
\begin{equation}\label{eq3}
\varphi(\mu(y, x),\alpha (z))-\varphi(\alpha(y),
\mu(x,z))=\mu(\alpha(y), \varphi(x,z))-\mu(\varphi(y, x),\alpha
(z)).
\end{equation}
By adding the equations \ref{eq1}, \ref{eq2} and \ref{eq3} and
considering the fact that $\varphi$ is skewsymmetric and $\mu$ is
commutative one has
\begin{equation*}
2\ \varphi(\alpha (x),\mu(y,z))=2\ \mu(\alpha(y),\varphi(x,z))+2\
\mu(\alpha(z),\varphi(x,y)).
\end{equation*}
\end{proof}
Let ${\A_t}=(V,\mu_t,\alpha_t)$ be a deformation
of the commutative Hom-associative algebra
${\A_0}=(V,\mu_0,\alpha_0)$. Assume that
\[
\mu_{t}(x, y) =\mu_{0}(x, y)+\mu_{1}(x, y)t
+\mu_{2}(x, y)t^{2}+\ldots,\quad \forall x,y\in V
\]
and
\[
\alpha_{t}(x) =\alpha_{0}(x)+\alpha_{1}(x)t
+\alpha_{2}(x,)t^{2}+\ldots,\quad \forall x\in V.
\]
 Then
$$
\frac{\mu_t (x,y)-\mu_t (y,x)}{t}=\mu_1 (x,y)-\mu_1 (y,x)+t
\sum_{i\geq 2}(\mu_{i}(x, y)-\mu_{i}(y, x))t^{i-1}.
$$
Hence, if $t$ goes to zero then $ \frac{\mu_t (x,y)-\mu_t (y,x)}{t}$
goes to $\{x ,y \}:=\mu_1 (x,y)-\mu_1 (y,x)$. The previous bracket
will define a structure of Poisson algebra over the commutative
algebra $\A_0$.

\begin{lem}\label{Compatibility}
Let  ${\A_0}=(V,\mu_0,\alpha_0)$ be a commutative Hom-associative
algebra and ${\A_t}=(V,\mu_t,\alpha_t)$ be a deformation of $\A_0$.
Then
$$\circlearrowleft_{x,y,z}{\delta^2_{Hom}\mu_{2}(x,y,z)}=
\circlearrowleft_{x,y,z}\mu_2\circ_{\alpha_0}\mu_0(x,y,z).$$
\end{lem}
\begin{proof}
\begin{align*}
\circlearrowleft_{x,y,z}{\delta^2_{Hom}\mu_{2}(x,y,z)}&=
\mu_2(\alpha_0(x), \mu_0(y,z))-\mu_2(\mu_0(x,y),\alpha_0(z)
\\ & + \mu_0(\alpha_0(x), \mu_2(y,z))-\mu_0(\mu_2(x,y),\alpha_0(z))\\
& +\mu_2(\alpha_0(y), \mu_0(z,x))-\mu_2(\mu_0(y,z),\alpha_0(x)
\\ & + \mu_0(\alpha_0(y), \mu_2(z,x))-\mu_0(\mu_2(y,z),\alpha_0(x))
\\
& +\mu_2(\alpha_0(z), \mu_0(x,y))-\mu_2(\mu_0(z,x),\alpha_0(y)
\\ &  +\mu_0(\alpha_0(z), \mu_2(x,y))-\mu_0(\mu_2(z,x),\alpha_0(y))
\\ &=  \circlearrowleft_{x,y,z}\mu_2\circ_{\alpha_0}\mu_0(x,y,z).
\end{align*}
\end{proof}

\begin{lem}
Let  ${\A_0}=(V,\mu_0,\alpha_0)$ be a commutative Hom-associative
algebra and ${\A_t}=(V,\mu_t,\alpha_0)$ be a deformation of $\A_0$.
Then
$$\circlearrowleft_{x,y,z}{\delta^2_{Hom}\mu_{2}(x,y,z)}-\circlearrowleft_{x,z,y}{\delta^2_{Hom}\mu_{2}(x,z,y)}=
0.$$
\end{lem}
\begin{proof}
\begin{align*}
\circlearrowleft_{x,y,z}{\delta^2_{Hom}\mu_{2}(x,y,z)}-\circlearrowleft_{x,z,y}{\delta^2_{Hom}\mu_{2}(x,z,y)}
&= \circlearrowleft_{x,y,z}\mu_2\circ_{\alpha_0}\mu_0(x,y,z)-
\circlearrowleft_{x,z,y}\mu_2\circ_{\alpha_0}\mu_0(x,z,y)\\ & =
\mu_2(\alpha_0(x), \mu_0(y,z))-\mu_2(\mu_0(x,y),\alpha_0(z))
\\ & + \mu_0(\alpha_0(y), \mu_2(z,x))-\mu_0(\mu_2(y,z),\alpha_0(x))
\\
& +\mu_2(\alpha_0(z), \mu_0(x,y))-\mu_2(\mu_0(z,x),\alpha_0(y))\\&
-\mu_2(\alpha_0(x), \mu_0(z,y))+\mu_2(\mu_0(x,z),\alpha_0(y))
\\ & - \mu_2(\alpha_0(z), \mu_0(y,x))+\mu_2(\mu_0(z,y),\alpha_0(x))
\\
& -\mu_2(\alpha_0(y), \mu_0(x,z))-\mu_2(\mu_0(y,x),\alpha_0(z))\\ &
=0.
\end{align*}
\end{proof}

\begin{lem}
Let  ${\A_0}=(V,\mu_0,\alpha_0)$ be a commutative Hom-associative
algebra and $\alpha$ any linear map of $Hom(V,V)$. Then
$$\circlearrowleft_{x,y,z}{\mu_{0}\circ_{\alpha}\mu_{0} (x,y,z)}=0.$$
\end{lem}
\begin{proof}
\begin{align*}
\circlearrowleft_{x,y,z}{\mu_{0}\circ_{\alpha}\mu_{0} (x,y,z)} &=
\mu_0(\alpha(x), \mu_0(y,z))-\mu_0(\mu_0(x,y),\alpha(z))
\\ & + \mu_0(\alpha(y), \mu_0(z,x))-\mu_0(\mu_0(y,z),\alpha(x))
\\
& +\mu_0(\alpha(z), \mu_0(x,y))-\mu_(\mu_0(z,x),\alpha(y))\\ &
=0.
\end{align*}
\end{proof}

\begin{thm}
Let  ${\A_0}=(V,\mu_0,\alpha_0)$ be a commutative Hom-associative
algebra and ${\A_t}=(V,\mu_t,\alpha_t)$ be a deformation of $\A_0$.
Consider the bracket  defined for $x,y\in V$ by $\{x ,y \}=\mu_1
(x,y)-\mu_1 (y,x)$ where $\mu_1$ is the first order element of the
deformation $\mu_t$.

Then $(V,\mu_0,\{ , \},\alpha_0)$ is a Hom-Poisson algebra.
\end{thm}
\begin{proof}
The bracket is skewsymmetric by definition. The compatibility
condition follows from  Lemma \ref{Compatibility}. Let us prove that
the Hom-Jacobi condition is satisfied by the bracket. One has
\begin{eqnarray*}
\circlearrowleft_{x,y,z}{\{\alpha_0(x),\{y,z\}\}}&=&
\circlearrowleft_{x,y,z}(\mu_1(\alpha_0(x),
\mu_1(y,z))-\mu_1(\alpha_0(x), \mu_1(z,y))+
\\ & &
-\mu_1(\mu_1(y, z),\alpha (x))+\mu_1(\mu_1(z,
y),\alpha (x))) \\ & =& \mu_1(\alpha_0(x),
\mu_1(y,z))-\mu_1(\alpha_0(x), \mu_1(z,y))+
\\ & &
-\mu_1(\mu_1(y, z),\alpha (x))+\mu_1(\mu_1(z,
y),\alpha (x))
\\ & & +\mu_1(\alpha_0(y),
\mu_1(z,x))-\mu_1(\alpha_0(y), \mu_1(x,z))+
\\ & &
-\mu_1(\mu_1(z,x ),\alpha (y))+\mu_1(\mu_1(x,
z),\alpha (y))
\\ & & +\mu_1(\alpha_0(z),
\mu_1(x,y))-\mu_1(\alpha_0(z), \mu_1(y,x))+
\\ & &
-\mu_1(\mu_1(x,y ),\alpha (z))+\mu_1(\mu_1(y,
x),\alpha (z))
\\ &=&\circlearrowleft_{x,y,z}{\mu_1
\circ_{\alpha_0} \mu_1(x,y,z)}-\circlearrowleft_{x,z,y}{\mu_1
\circ_{\alpha_0} \mu_1(x,z,y)}.
\end{eqnarray*}
The deformation equation \eqref{DefEqu}  implies for $s=2$ that
$$\mu_{1}\circ_{\alpha_0} \mu_{1}=-\mu_{2}\circ_{\alpha_0} \mu_{0}-\mu_{0}\circ_{\alpha_0}
\mu_{2}-\mu_{1}\circ_{\alpha_1} \mu_{0}-\mu_{0}\circ_{\alpha_1}
\mu_{1}-\mu_{0}\circ_{\alpha_2} \mu_{0}$$ which is equivalent to
$$\mu_{1}\circ_{\alpha_0} \mu_{1}=-\delta^2_{Hom}\mu_{2}-\mu_{0}\circ_{\alpha_1} \mu_{1}-\mu_{1}\circ_{\alpha_1}
\mu_{0}-\mu_{0}\circ_{\alpha_2} \mu_{0}$$ Then using the previous
Lemmas
\begin{align*}\circlearrowleft_{x,y,z}{\{\alpha(x),\{y,z\}\}}&=
\circlearrowleft_{x,y,z}{(-\delta^2_{Hom}\mu_{2}(x,y,z)-\mu_{0}\circ_{\alpha_1}
\mu_{1}(x,y,z)-\mu_{1}\circ_{\alpha_1}
\mu_{0}(x,y,z)-\mu_{0}\circ_{\alpha_2} \mu_{0}(x,y,z))}\\&
+\circlearrowleft_{x,z,y}{(\delta^2_{Hom}\mu_{2}(x,z,y)+\mu_{0}\circ_{\alpha_1}
\mu_{1}(x,z,y)+\mu_{1}\circ_{\alpha_1}
\mu_{0}(x,z,y)+\mu_{0}\circ_{\alpha_2} \mu_{0}(x,z,y))}\\ &  =0
\end{align*}
\end{proof}

\section{On Cohomolgy of Hom-Lie algebras}
Now we introduce elements of cohomology of Hom-Lie algebras in connection to their infinitesimal deformations.

\subsection{First and second coboundary operators}

 Let ${\G }=(V,[\cdot, \cdot ],\alpha)$ be a Hom-Lie algebra on
a $\K$-vector space $V.$

The set of $p$-cochains on $V$ is the set of $p$-linear alternating
maps
\begin{equation*}
\mathcal{C}^{p}({\G,\G})=\{\varphi:V^{\wedge p} =\underset{p\text{
times}}{\underbrace{V\wedge V\wedge ...\wedge V}}\longrightarrow V\}
\end{equation*}
In the following, we define $1$-Hom-cochains and $2$-Hom-cochains.
\begin{definition}
A $1$-\emph{Hom-cochain} is a map $f$, where $f\in
\mathcal{C}^{1}({\G,\G})$  satisfying
\begin{equation}
f\circ\alpha=\alpha\circ f
\end{equation}
We denote by $Hom\mathcal{C}^{1}({\G,\G})$ the set of all
$1$-\emph{Hom-cochain} of $\G$.
\end{definition}

\begin{definition}
A $2$-\emph{Hom-cochain} is a pair $(\varphi, \tau)$, where
$\varphi\in \mathcal{C}^{2}({\G,\G})$ is a $2$-linear alternating
map and $\tau$ is a linear map satisfying
\begin{equation}
\circlearrowleft_{x,y,z}{[\tau(x),[y,z]]}=0
\end{equation}
We denote by $Hom\mathcal{C}^{2}({\G,\G})$ the set of all
$2$-\emph{Hom-cochains} of $\G$.
\end{definition}

The $1$-coboundary and $2$-coboundary operators for Hom-Lie algebras
are defined as follows

\begin{definition}
We call  $1$-\emph{coboundary} operator of Hom-Lie algebra ${\G}$
the map
\begin{equation*}
\delta_{HL} ^{1}:  \mathcal{C} ^{1}({\G,\G})  \longrightarrow
\mathcal{C}^{2}({\G ,\G }), \quad f  \longmapsto \delta_{HL}^{1}f
\end{equation*}
defined by
\begin{equation*}
\delta_{HL} ^{1}f(x, y) =f([x,y])-[f(x), y]- [x, f(y)]
\end{equation*}

\end{definition}

\begin{definition}
We call a $2$-\emph{coboundary} operator of Hom-Lie algebra ${\G}$
the map
\begin{equation*}
\delta_{HL}^{2}:  \mathcal{C}^{2}({\G ,\G }) \longrightarrow
\mathcal{C}^{3}({\G ,\G }), \quad
 \varphi  \longmapsto  \delta_{HL} ^{2}\varphi
\end{equation*}
defined by
\begin{align*}
\delta_{HL}^{2}\varphi (x, y, z) &=
\circlearrowleft_{x,y,z}{\varphi(\alpha(x),[y,z])+[\alpha(x),\varphi(y,z)]}
\end{align*}
\end{definition}

The cohomology spaces relative to these coboundary operators are

\begin{definition}
The space of $1$-cocycles of ${\A}$ is
\begin{equation*}
\mathit{H_{HL}^{1}}({\G,\G})=\{f\in
Hom\mathcal{C}^{1}({\G ,\G
})\,:\,\delta_{HL}^{1}f=0\ \}=\{f\in
\mathcal{C}^{1}({\G ,\G
})\,:\,\delta_{HL}^{1}f=0\ \text{ and }
f\circ\alpha=\alpha\circ f \}
\end{equation*}
The space of $2$-coboundaries of ${\A }$ is
\begin{equation*}
\mathit{B_{HL}^{2}}({\G ,\G }) =\{(\varphi,\tau)\in
Hom\mathcal{C}^{2}({\G ,\G }) \ :\, \varphi =\delta_{HL} ^{1}f,\
f\in Hom\mathcal{C} ^{1}({\G },{\G })\ \text{ and }
\circlearrowleft_{x,y,z}{[\tau(x),[y,z]]}=0\ \forall x,y,z\in V \}
\end{equation*}
The space of $2$-cocycles of ${\G}$ is
\begin{equation*}
\mathit{Z_{HL}^{2}}({\G ,\G }) =\{(\varphi,\tau)\in
Hom\mathcal{C}^{2}({\G ,\G })\,\,:\,\delta_{HL} ^{2}\varphi=0 \
\text{ and } \circlearrowleft_{x,y,z}{[\tau(x),[y,z]]}=0\ \forall
x,y,z\in V \}
\end{equation*}
\end{definition}
\begin{proposition}
$\delta_{HL} ^{2}( \delta_{HL} ^{1})=0$
\end{proposition}
\begin{proof}
Let
\begin{equation*}
\delta_{HL} ^{1}f(x, y) =f([x,y])-[f(x), y]- [x,
f(y)]
\end{equation*}
One has
\begin{align*}
\delta_{HL}^{2}(\delta_{HL} ^{1} f)(x, y, z) &=f([\alpha
(x),[y,z]])-[f(\alpha (x)),[y,z]]-[\alpha (x),f([y,z])]\\
& +[\alpha
(x),f([y,z])])-[\alpha (x)),[f(y),z]]-[\alpha (x),[y,f(z)]]\\
& f([\alpha
(y),[z,x]])-[f(\alpha (y)),[z,x]]-[\alpha (y),f([z,x])]\\
& +[\alpha
(y),f([z,x])])-[\alpha (y)),[f(z),x]]-[\alpha (y),[z,f(x)]]\\
& f([\alpha
(z),[x,y]])-[f(\alpha (z)),[x,y]]-[\alpha (z),f([x,y])]\\
& +[\alpha (z),f([x,y])]-[\alpha (z)),[f(x),y]]-[\alpha
(z),[x,f(y)]].
\end{align*}
Since $\alpha$ and $f$ commute and the bracket satisfies the
Hom-Jacobi identity, then
\begin{eqnarray*}
\delta_{HL}^{2}(\delta_{HL} ^{1} f)(x, y, z) &=&
f(\circlearrowleft_{x,y,z}{[\alpha
(x),[y,z]]})-\circlearrowleft_{f(x),y,z}{[\alpha (f( x)),[y,z]]} \\
 & & -\circlearrowleft_{x,f(y),z}{[\alpha
(x),[f(y),z]]}-\circlearrowleft_{x,y,f(z)}{[\alpha (x),[y,f(z)]]}=0.
\end{eqnarray*}

\end{proof}
\begin{remark}
One has $\mathit{B_{HL}^{2}}({\G ,\G }){\subset
}\mathit{Z_{HL}^{2}}({\G ,\G })$, because $\delta_{HL} ^{2}\circ
\delta_{HL} ^{1}=0$. Note also that $\mathit{Z_{HL}^{1}}({\G ,\G })$
gives the space of derivations of a Hom-Lie algebra ${\G }$, denoted
$Der_{HL}{(\G )}$.
\end{remark}

\begin{definition}
We call the $2 ^{th}$ cohomology group of the
Hom-Lie algebra ${\A }$ the quotient
\begin{equation*}
\mathit{H_{HL}^{2}}(\G ,\G )
=\frac{\mathit{Z_{HL}^{2}}(\G
,\G)}{\mathit{B_{HL}^{2}}(\G ,\G )}.
\end{equation*}
\end{definition}
\subsection{Deformations of Hom-Lie algebras in terms of cohomology}
The results in this section are similar to those obtained for
Hom-associative algebra in Section 4.

Let ${\G }=(V,[\cdot, \cdot ],\alpha)$ be a Hom-Lie algebra on a
$\K$-vector space $V$ and ${\G_t }=(V,[\cdot,\cdot]_t,\alpha_t)$
where $ [x,y]_{t} =\sum_{i\geq 0}[x,y]_{i}t^{i}$ and $ \alpha_{t}(x)
=\sum_{i\geq 0}\alpha_{i}(x)t^{i}$, $x,y\in V$, such that
$[\cdot,\cdot]_i$ are bilinear alternating maps with
$[\cdot,\cdot]_0$ being the bracket of $\G$ and $\alpha_0=\alpha$.
We assume in this Section that the deformation satisfies
$\circlearrowleft_{x,y,z}{[\alpha_1(x),[y,z]_0]_0}=0$.
\begin{proposition}
The first order  term $([\cdot,\cdot]_1,\alpha_1)$ of the
deformation ${\G_t }=(V,[\cdot,\cdot]_t,\alpha_t)$ of ${\G
}=(V,[\cdot, \cdot ],\alpha)$ satisfying
$\circlearrowleft_{x,y,z}{[\alpha_1(x),[y,z]_0]_0}=0$ is a
2-Hom-cocycle of  the Hom-Lie algebra $\G$.
\end{proposition}
\begin{proof}
The deformation equation \ref{equa2} implies
$$\sum_{i+j+k=s}\circlearrowleft_{x,y,z}{[\alpha_j(x),[y,z]_k]_i}=0, \text{ for } s>0
$$
For $s=1$, one has
$$\circlearrowleft_{x,y,z}{[\alpha_0(x),[y,z]_0]_1+[\alpha_1(x),[y,z]_0]_0+[\alpha_0(x),[y,z]_1]_0}=0.
$$
Since $\circlearrowleft_{x,y,z}{[\alpha_1(x),[y,z]_0]_0}=0$ then for
$\psi=[\cdot,\cdot]_1$, $\delta_{HL} ^{2}\psi=0$.
\end{proof}

\section{Examples of Hom-Lie deformations of   $\sll$}
In this section, we construct all the twistings so that the  $\sll$ brackets determine a Hom-Lie algebra,   provide families of Hom-Lie algebras deforming  the $\sll$ Lie algebra and show that $q$-deformed Witt
algebras may be viewed as a Hom-deformation.
\subsection{Infinitesimal Hom-Lie deformations of $\sll$} In this section, we deform the $\sll$ Lie algebra as a Hom-Lie algebra.  The following proposition gives  all the twistings
so that the brackets $[X_1,X_2]=2 X_2, \
[X_1,X_3]=-2 X_3, \ [X_2,X_3]=X_1$ determine a
three dimensional Hom-Lie algebra, generalizing
$\sll$.

\begin{proposition}
Let $V$ be a three dimensional $\K$-linear space and let $(x_1 ,x_2,
x_3)$ be its basis. Any Hom-Lie algebra with the following brackets
$$[x_1,x_2]=2 x_2, \ [x_1,x_3]=-2 x_3, \ [x_2,x_3]=x_1
$$ is given by linear $\alpha$ defined, with respect to the previous basis, by a matrix of the form
 $\left(
  \begin{array}{ccc}
  a & d& c \\
  2c & b & f \\
    2d & e & b\\
  \end{array}
\right) $ where $a,c,d,e,f$ are arbitrary parameters in $\K$

\end{proposition}
\begin{remark}
The Hom-Lie algebras given by the previous theorem are deformations of $\sll$ viewed as a Hom-Lie algebra where $\alpha_0$ is the identity matrix.
\end{remark}

In the following we provide infinitesimal Hom-Lie deformations of $\sll$. We construct $3$-dimensional Hom-Lie algebras defined by the bracket $[~,~]_t=[~,~]_0+t\ [~,~]_1$ and homomorphism $\alpha_t=\alpha_0+t\ \alpha_1$, where $\alpha_0$ is the identity map,   satisfying Hom-Jacobi condition and such that $( [~,~]_1,\alpha_1)$ is a 2-Hom-cocycle for the  cohomology Hom-Lie algebras defined above.

The following pairs $( [~,~]_1,\alpha_1)$ define a 2-Hom-cocyle and  thus corresponding infinitesimal deformations of $\sll$.

\begin{enumerate}
\item $
\begin{array}{cc}
\begin{array}{ccc}
 [ x_1, x_2 ] _1&= &-a_1 x_2 +x_3 \\ {}
 [x_1, x_3 ]_1&=& a_2 x_2+a_1 x_3, \\ {}
 [ x_2,x_3 ] _1& = & a_3 x_1,
 \end{array}
 & \quad
 \alpha_1=\left(
  \begin{array}{ccc}
  b_1 & 0& 0 \\
  0 & b_2 & -b_3 a_2 \\
    0 & b_3 & b_2
  \end{array}
\right)
\end{array}
$

where $ a_1,a_2,a_3,b_1,b_2$ are parameters in $\K$. \vspace{0.5cm}

\item $
\begin{array}{cc}
\begin{array}{ccc}
 [ x_1, x_2 ] _1&= &-2 a_1 x_2  \\ {}
 [x_1, x_3 ]_1&=& a_2 x_2+2 a_1 x_3, \\ {}
 [ x_2,x_3 ] _1& = & -a_1 x_1,
 \end{array}
 & \quad
 \alpha_1=\left(
  \begin{array}{ccc}
  b_1 & 0& b_3 \\
2 b_3 & b_2 & b_4 \\
    0 & 0 & b_2
  \end{array}
\right)
\end{array}
$

where $ a_1,a_2,b_1,b_2,b_3,b_4$ are parameters in $\K$.
\vspace{0.5cm}

\item $
\begin{array}{cc}
\begin{array}{ccc}
 [ x_1, x_2 ] _1&= &- a_1 x_2+a_2 x_3  \\ {}
 [x_1, x_3 ]_1&=& a_3 x_1+a_4 x_2+ a_1 x_3, \\ {}
 [ x_2,x_3 ] _1& = & a_5 x_1-a_3 x_2,
 \end{array}
 & \quad
 \alpha_1=\left(
  \begin{array}{ccc}
  b & 0& 0 \\
0 & b & 0 \\
    0 & 0 & b
  \end{array}
\right)
\end{array}
$

where $ a_1,a_2,a_3,a_4,a_5,b$ are parameters in $\K$ with $a_3\neq
0$. \vspace{0.5cm}
\end{enumerate}

\begin{remark}
The above three examples are actually Lie algebras for all values of parameters $a_j$ and $b_j$
since $b_j's$ can be chosen arbitrary and in particular so that $\alpha_1 = id$ and hence
$\alpha_0 + t \alpha_1 = (1+t) id $ corresponding to the case of Lie algebra.
We have also computed with a  computer algebra  system "Mathematica" many other examples of infinitesimal formal   Hom-Lie  deformations of $\sll$
with additional  restriction that $(V,[~,~]_0,\alpha_1)$ is a Hom-Lie algebra. Remarkably,  it turns out that in all these examples the Hom-Lie algebras are actually Lie algebras. We conjecture that this is always the  case for such Hom-Lie infinitesimal formal deformations of $\sll$.
\end{remark}

Now, we will give examples of Hom-Lie infinitesimal formal deformations of $\sll$ which are not Lie algebras.  We consider the  $3$-dimensional Hom-Lie algebras with the bracket $[~,~]_t$ and linear map $\alpha_t$  defined  as follows
$$
\begin{array}{cc}
\begin{array}{ccc}
 [ x_1, x_2 ] _t&= &a_1 t x_1+(2- a_2 t) x_2  \\ {}
 [x_1, x_3 ]_t&=& a_3 t x_1+a_4 t x_2+(-2+ a_2 t) x_3, \\ {}
 [ x_2,x_3 ] _t& = & (1-\frac{a_2}{2} t) x_1,
 \end{array}
 & \ \
 \alpha_t=\left(
  \begin{array}{ccc}
 1+ b_1 t & \frac{a_1}{2}t& \frac{b_2-a_3}{2}t \\
 b_2 t &1- \frac{a_2 }{2}t&-\frac{a_4}{2}t \\
    0 & 0 &1- \frac{a_2 }{2}
  \end{array}
\right)
\end{array}
$$

where $ a_1,a_2,a_3,a_4,b_1,b_2$ are parameters in $\K$. This Hom-Lie algebra becomes a Lie algebra for all $t$ if and only if $a_1=0$ and $a_3=0$, as follows from
$$[x_1,[x_2,x_3]]+[x_3,[x_1,x_2]]+[x_2,[x_3,x_1]]=
(2a_3t -(a_2a_3 +a_1a_4)t^2)x_2+(2a_1t-a_1a_2 t^2)x_3.
$$

\subsection{$q$-deformed Witt algebras}

 Let $\A$ be the unique
factorization domain $\mathbb{K}[z,z^{-1}]$, the
Laurent polynomials in $t$ over the field
$\mathbb{K}$. Then the space ${\mathcal
D}_\sigma(\A)$ can be generated by a single
element $D$ as a left $\A$-module, that is,
${\mathcal D}_\sigma(\A)=\A\cdot D$ (Theorem 1 in
\cite{HLS}). When $\sigma(z)=qz$ with $q\neq0$
and $q\neq 1$, one can take $D$ as $z$ times the
Jackson $q$-derivative
$$D=\frac{id-\sigma}{1-q}\quad :\quad
f(t)\mapsto\frac{f(z)-f(qz)}{1-q}.$$
The $\mathbb{K}$-linear space
    ${\mathcal
D}_\sigma (\A) = \bigoplus_{n\in \mathbb{Z}}
\mathbb{K}\cdot x_n,$ with
    $x_n=-z^nD$ can be equipped with the skew-symmetric
    bracket $[ \cdot,\cdot ]_\sigma$
    defined on generators  as
    \begin{equation} \label{qWittrel} [x_n,x_m]=
    q^nx_{n}x_m-q^mx_{m}x_n =[x_n,x_m]
    =(\{n\}_q-\{m\}_q)x_{n+m},
    \end{equation}
    where $\{n\}_q=(q^n-1)/(q-1)$ for $q\neq 1$
    and $\{n\}_1=n.$
    This bracket is skew-symmetric and
    satisfies the $\sigma$-deformed Jacobi-identity
    \begin{equation} \label{qWittJacobi}
        (q^n+1)[ x_n, [ x_l, x_m]
        ]+ (q^l+1)[ x_l,  [ x_m, x_n
]_\sigma]
    +(q^m+1)[ x_m, [ x_n,x_l]
    ]=0.
    \end{equation}
    We have a Hom-Lie algebra
    $(V,[\cdot,\cdot],\alpha) =
    (\bigoplus_{i\in \mathbb{Z}} \K
    x_n, [\cdot,\cdot],\alpha)$
    with the bilinear bracket defined on generators as
$[x_n,x_m] = (\{n\}_q-\{m\}_q)x_{n+m}$ and
   the linear twisting map $\alpha:V\rightarrow V$
    acting on generators
    as $\alpha (x_n) =  (q^n+1) x_n$.
    Obviously this Hom-Lie algebra can be viewed
    as a $q$-deformed Witt algebra in the sense that
    for $q = 1$ indeed one recovers the bracket and
    the commutation
    relations for generators of the Witt algebra.
    The definition of its generators using first
    order differential operators is recovered if one assumes that $D= t\frac{d}{dt}$
    for $q=1$ as one would expect from
    passing to a limit
    in the definition of the operator $D$.

    It can be also shown that there is a central
    extension $Vir_q$ of this deformation in the category of hom-Lie
    algebras \cite{HLS}, therefore being a natural $q$-deformation
    of the Virasoro algebra. The algebra $Vir_q$ is
    spanned by elements
    $\{x_n\,|\,n\in\mathbb{Z}\}\cup\{{\bf c}\}$ where
    ${\bf c}$ is central with respect to Hom-Lie bracket,
    i.e.,
    $[ Vir_q,{\bf c}]=
    [{\bf c}, Vir_q]=0$. The
    bracket of $x_n, x_m$ is computed according to
    $$[ x_n,x_m]=(\{n\}_q-\{m\}_q)x_{m+n}+
    \delta_{n+m,0} \frac{q^{-n}}{6(1+q^n)}\{n-1\}_q\{n\}_q\{n+1\}_q{\bf c}.$$
    Note that when $q=1$ we retain the
         classical Virasoro algebra
         $$[ x_n,x_m]=(n-m)x_{m+n}+
         \delta_{n+m,0}
         \frac{1}{12}\{n-1\}\{n\}\{n+1\} {\bf c}.$$
         from conformal field and string theories.
         Note also that when specializing ${\bf c}$ to
         zero, or equivalently rescaling ${\bf c}$
         by extra parameter and then letting the
         parameter degenerate to zero,
         one recovers the $q$-deformed Witt
         algebra. Because of this the Witt
         algebra is called also,
         primarily in the physics literature,
a centerless Virasoro algebra. In a similar way
the $q$-deformed Witt algebra could be called a
centerless $q$-deformed Virasoro algebra but of
course with the word "central" used in terms of
the Hom-Lie algebra bracket.

Next, we will link the $q$-deformed Witt algebras
to the framework of deformation theory of Hom-Lie
algebras as we have developed. To this end, let $q=1+t$. For $n\in
\mathbb{Z}_{\geq 0} = \{0,1,2,\dots\}$ we have

\begin{eqnarray*}
\{n\}_q &=& \sum_{j=0}^{n-1} q^j =
\sum_{j=0}^{n-1} (1+t)^j = \sum_{j=0}^{n-1}
\sum_{k=0}^{j} \binom{j}{k} t^k =
\sum_{k=0}^{n-1}
(\sum_{j=k}^{n-1} \binom{j}{k})t^k \\
q^n &=& (1+t)^n= \sum_{k=0}^{n} \binom{n}{k} t^k,
\end{eqnarray*}
and thus $\alpha_t (x_n) =  (q^n + 1) x_n = (1+
\sum_{k=0}^{n} \binom{n}{k} t^k) x_n = (1+
\sum_{k=0}^{n} \binom{n}{k} x_n t^k). $ So, in
the decomposition
\begin{equation*} \alpha_{t}(x)
=\sum_{k\geq 0}\alpha_{k}(x)t^{k},\quad
\text{where} \quad \alpha_{k}\in Hom(V,V)
\end{equation*}
we have $\alpha_{0}(x_n) = 2x_n$ and
$\alpha_{k}(x_n) = \binom{n}{k} x_n$ for $k>0$
and $n\in \mathbb{Z}$. In particular,
$\alpha_{1}(x_n) = \binom{n}{1} x_n = n x_n$

For $n, m \in \mathbb{Z}_{\geq 0}$, the relations
\eqref{qWittrel} defining the Hom-Lie bracket on
generators can be rewritten in terms of $t$ as
follows: \begin{align*} [ x_n,x_m]_t &=
\sum_{k=0}^{\max(n,m)-1}
\left(\left(\sum_{j=k}^{n-1} \binom{j}{k} -
\sum_{j=k}^{m-1}
\binom{j}{k}\right)x_{n+m}\right) t^k \\
&= \sum_{k\in \mathbb{Z}_{\geq
0}}[x_n,x_m]_{k}t^{k},
\end{align*}
 where
$$[x_n,x_m]_{k} =\left(\sum_{j=k}^{n-1} \binom{j}{k} -
\sum_{j=k}^{m-1} \binom{j}{k}\right)x_{n+m}
\text{ \ for \ } k\in \mathbb{Z}_{\geq 0}.$$ In
particular, for $k=0$ we have
\begin{equation*}
[x_n,x_m]_0 = \left(\sum_{j=0}^{n-1} \binom{j}{0}
- \sum_{j=0}^{m-1} \binom{j}{0}\right)x_{n+m} =
(n-m)x_{n+m}
\end{equation*}
meaning that the Witt algebra is exactly present
in the zero degree term (origin) of the
deformation. In the first order term, $k=1$, we
get \begin{align*} [x_n,x_m]_1 &=
\left(\sum_{j=1}^{n-1} \binom{j}{1} -
\sum_{j=1}^{m-1} \binom{j}{1}\right)x_{n+m} \\
&=
\left(\frac{n(n-1)}{2}-\frac{m(m-1)}{2}\right)x_{n+m}
= \frac{(n-m)(n+m-1)}{2} x_{n+m}.\end{align*}

\begin{proposition}
Consider the Witt algebra $W_{\geq
0}=\bigoplus_{n\in \mathbb{Z}_{\geq 0}} \K x_n$
defined by the brackets
$$[x_n,x_m]_0=(n-m) x_{n+m}.$$

The one parameter families given by

\begin{align*} [ x_n,x_m]_t = \sum_{k\in
\mathbb{Z}_{\geq 0}}[x_n,x_m]_{k}t^{k} \quad
\text{ and } \quad \alpha_{t}(x) =\sum_{k\in
\mathbb{Z}_{\geq 0}}\alpha_{k}(x)t^{k}
\end{align*}
where
\begin{align*}[x_n,x_m]_{k} &=\left(\sum_{j=k}^{n-1}
\binom{j}{k} - \sum_{j=k}^{m-1}
\binom{j}{k}\right)x_{n+m} \text{ \ for \ } k\in
\mathbb{Z}_{\geq 0},\\ \alpha_{0}(x_n) &= 2x_n,
\quad \alpha_{k}(x_n) = \binom{n}{k} x_n \text{
for } k>0, n\in \mathbb{Z} \end{align*} define a
Hom-Lie algebra deformation of the Witt algebra
$W_{\geq 0}$.

\end{proposition}
\begin{remark}
Consider $W_{\geq 0}=\bigoplus_{n\in
\mathbb{Z}_{\geq 0}} \K x_n$. Let
$[\cdot,\cdot]_1$ be an skew-symmetric bilinear
map on $W_{\geq 0}$ and $\alpha_1$ be a linear
map on $W_{\geq 0}$, defined by
$$
[x_n,x_m]_1=\frac{(n-m)(n+m-1)}{2}x_{n+m}
$$
and
$$\alpha_1(x_n)=n x_{n}.$$
Then the pair $([\cdot,\cdot]_1,\alpha_1)$ is not a
2-Hom-cocyle (with respect to the previous cohomology) of  the Witt algebra $W_{\geq 0}$
considered as a Hom-Lie algebra with the linear
map defined by $\alpha_0(x_n)=2 x_{n}$.

One has
\begin{eqnarray*}
\circlearrowleft_{p,q,r}{
[\alpha_0(x_p),[x_q,x_r]_0]_1
+[\alpha_1(x_p),[x_q,x_r]_0]_0+
[\alpha_0(x_p),[x_q,x_r]_1]_0} \\=
\circlearrowleft_{p,q,r}
{[2x_p,(q-r)x_{q+r}]_1+[p x_p,(q-r) x_{q+r}]_0+
[2x_p,\frac{(q-r)(q+r-1)}{2}x_{q+r}]_0} \\=
\circlearrowleft_{p,q,r}{2(q-r)(p-q-r)(p+q+r-1)
x_{p+q+r}}=0
\end{eqnarray*}
but
$$\circlearrowleft_{p,q,r}{
[\alpha_0(x_p),[x_q,x_r]_0]_1+
[\alpha_0(x_p),[x_q,x_r]_1]_0} \neq 0$$
\end{remark}

\end{document}